\documentstyle[]{article}

\def\qmod#1#2{{\hbox{}^{\displaystyle{#1}}}\!\big/\!\hbox{}_{
\displaystyle{#2}}}
\def\catqmod#1#2{{\hbox{}^{\displaystyle{#1}}}\!\big/\hskip-4pt{\big/}\!\hbox{}_{
\displaystyle{#2}}}

\newfam\msbfam

\font\twelmsb=msbm10 at 12pt
\font\tenmsb=msbm10

\font\sevenmsb=msbm10 at 7pt
\font\fivemsb=msbm10 at 5pt

\textfont\msbfam=\tenmsb
\scriptfont\msbfam=\sevenmsb
\scriptscriptfont\msbfam=\fivemsb

\def\Bbb{\fam\msbfam\twelmsb}

\def\A{{\Bbb A}}

\def\C{{\Bbb C}}

\def\F{{\Bbb F}}
\def\G{{\Bbb G}}

\def\P{{\Bbb P}}
\def\Q{{\Bbb Q}}
\def\R{{\Bbb R}}
\def\Z{{\Bbb Z}}

\def\qed {\hfill\vrule height6pt width6pt depth0pt \bigskip}

\def\map{\longrightarrow}
\def\textmap#1{\mathop{\vbox{\ialign{
                                ##\crcr
    ${\scriptstyle\hfil\;\;#1\;\;\hfil}$\crcr
    \noalign{\kern-1pt\nointerlineskip}
    \rightarrowfill\crcr}}\;}}

\def\textlmap#1{\mathop{\vbox{\ialign{
                                ##\crcr
    ${\scriptstyle\hfil\;\;#1\;\;\hfil}$\crcr
    \noalign{\kern-1pt\nointerlineskip}
    \leftarrowfill\crcr}}\;}}

\newfam\meuffam

\font\tenmeuf=eufm10
\font\sevenmeuf=eufm7
\font\fivemeuf=eufm5

\textfont\meuffam=\tenmeuf
\scriptfont\meuffam=\sevenmeuf
\scriptscriptfont\meuffam=\fivemeuf

\def\germ{\fam\meuffam\tenmeuf}

\def\g{{\germ g}}

\def\kg{{\germ k}}

\def\pg{{\germ p}}

\def\tg{{\germ t}}

\def\vg{{\germ v}}

\begin{document}
\def\Pr{{\rm Pr}}
\def\tr{{\rm Tr}}
\def\End{{\rm End}}
\def\Aut{{\rm Aut}}
\def\Spin{{\rm Spin}}
\def\U{{\rm U}}
\def\SU{{\rm SU}}
\def\SO{{\rm SO}}
\def\PU{{\rm PU}}
\def\GL{{\rm GL}}
\def\spin{{\rm spin}}
\def\u{{\rm u}}
\def\su{{\rm su}}
\def\so{{\rm so}}
\def\ub{\underbar}
\def\pu{{\rm pu}}
\def\Pic{{\rm Pic}}
\def\Iso{{\rm Iso}}
\def\NS{{\rm NS}}
\def\deg{{\rm deg}}
\def\Hom{{\rm Hom}}
\def\Aut{{\rm Aut}}
\def\h{{\germ h}}
\def\Herm{{\rm Herm}}
\def\Vol{{\rm Vol}}
\def\pf{{\bf Proof: }}
\def\id{{\rm id}}
\def\i{{\germ i}}
\def\im{{\rm im}}
\def\rk{{\rm rk}}
\def\ad{{\rm ad}}
\def\h{{\bf H}}
\def\coker{{\rm coker}}
\def\dv{\bar\partial}
\def\Ad{{\rm Ad}}
\def\RSU{\R SU}
\def\ad{{\rm ad}}
\def\dva{\bar\partial_A}
\def\da{\partial_A}
\def\p{\partial\bar\partial}
\def\sp{\Sigma^{+}}
\def\sm{\Sigma^{-}}
\def\spm{\Sigma^{\pm}}
\def\smp{\Sigma^{\mp}}
\def\oo{{\scriptstyle{\cal O}}}
\def\ooo{{\scriptscriptstyle{\cal O}}}
\def\sw{Seiberg-Witten }
\def\pa{\partial_A\bar\partial_A}
\def\Dr{{\raisebox{0.15ex}{$\not$}}{\hskip -1pt {D}}}
\def\gr{{\scriptscriptstyle|}\hskip -4pt{\g}}
\def\subsetint{{\  {\subset}\hskip -2.45mm{\raisebox{.28ex}
{$\scriptscriptstyle\subset$}}\ }}
\def\nr{\parallel}
\def\ra{\rightarrow}
\newtheorem{sz}{Satz}[section]
\newtheorem{thry}[sz]{Theorem}
\newtheorem{pr}[sz]{Proposition}
\newtheorem{re}[sz]{Remark}
\newtheorem{co}[sz]{Corollary}
\newtheorem{dt}[sz]{Definition}
\newtheorem{lm}[sz]{Lemma}
\newtheorem{cl}[sz]{Claim}

\title{Gauge theoretical  Gromov-Witten invariants and virtual fundamental classes}
\author{Ch. Okonek$^*$  
\and A. Teleman\thanks{Partially supported by: EAGER -- European Algebraic
Geometry Research Training Network, contract No HPRN-CT-2000-00099 (BBW
99.0030), and by SNF, nr. 2000-055290.98/1}}
\maketitle

\begin{abstract}\vspace{2mm}

This article is an expanded version of talks which the authors have given
in Oberwolfach, Bochum, and at the Fano Conference in Torino.
In these talks we explained the main results of our papers
"Gauge theoretical equivariant Gromov-Witten invariants and the full
Seiberg-Witten invariants of ruled surfaces" and "Comparing virtual
fundamental classes: Gauge theoretical Gromov-Witten invariants for toric
varieties".

We  have also included new results, e. g. the material concerning flag varieties,   Quot spaces
over $\P^1$, and the generalized quiver representations.

The common theme   is the construction of gauge theoretical
Gromov-Witten type invariants of arbitrary genus associated with certain
symplectic factorization problems with additional symmetry, and the
computation of these invariants in terms of complex geometric objects.

In chapter 1 we introduce the concept of a symplectic factorization problem
with additional symmetry (SFPAS), and illustrate it with several important
examples:
Grassmann manifolds, flag varieties, toric varieties, certain Quot spaces
over $\P^1$, and generalized quivers.

Chapter 2 introduces the gauge theoretical problem associated with a SFPAS,
i.e. the standard gauge theoretical set up consisting of a configuration
space, a partial differential equation (of vortex type), and a gauge group.
Our invariants are defined by evaluating canonical cohomology classes on
the virtual fundamental class of the moduli space of irreducible solutions
to the PDE when this is possible, e.g. when the moduli space is compact.

In chapter 3 we explain the complex geometric interpretation of these
moduli spaces in terms of (poly-)stable framed holomorphic objects over a
Riemann surface, provided the original SFPAS came from a nice K\"ahlerian
problem.
The point is that, in this case, the moment map of the SFPAS together with
the  Riemannian metric on the base surface gives rise to a naturally
associated stability concept for these framed holomorphic objects.
The "universal Kobayashi-Hitchin correspondence" shows then the existence
of an isomorphism $\iota:{\cal M}^*\ra{\cal M}^{\rm st}$ between the gauge theoretic moduli
spaces of irreducible solutions of the PDE and the complex geometric moduli space
of stable framed holomorphic objects.
However, in order to compute the gauge theoretical Gromov-Witten invariants
using this explicit desciption, one also needs to know that the
Kobayashi-Hitchin isomorphism $\iota:
{\cal M}^*\ra{\cal M}^{\rm st}$ identifies the virtual
fundamental classes of these moduli spaces. In general this is a very
difficult problem.
We show that this is true for the special SFPAS which yields the   toric
varieties, and we state a conjecture for the general situation. Roughly speaking this
conjecture asserts the following:
When the gauge theoretic problem is of Fredholm type, and the data for        ${\cal M}^{\rm
st}$ are algebraic, then ${\cal M}^{\rm st}$ admits a canonical perfect obstruction theory in
the sense of Behrend-Fantechi, and the Kobayashi-Hitchin isomorphism $\iota:
{\cal M}^*\ra{\cal M}^{\rm st}$ identifies the gauge theoretic and the algebraic virtual
fundamental classes.

Chapter 4 contains examples, explicit computations of our invariants in an
abelian case, as well as some interesting applications to Seiberg-Witten
invariants of ruled surfaces. Another nice application concerns an old
enumerative problem, namely the counting of maximal subbundles of a general
vector bundle over a Riemann surface.

\end{abstract}

\section{Symplectic factorization problems with additional symmetry}

A {\it symplectic factorization problem} (SFP) is a system $(F,\alpha,\mu)$,  where
$F=(F,\omega)$ is a symplectic manifold, $\alpha:K\times F\ra F$ is a symplectic action of a
compact Lie group $K$ on $F$ and $\mu$ is a moment map for this action.  The
{\it result} of a symplectic factorisation problem $(F,\alpha,\mu)$ is the quotient 
$$F_\mu:= \qmod{\mu^{-1}(0)}{K}\ .
$$
This quotient becomes a symplectic manifold (respectiveley orbifold) if $\mu$ is
a submersion in any point of $\mu^{-1}(0)$ and every point of this set   has
trivial (repectively finite) stabilizers with respect to the $K$-action.  In the first
case we will say that the SFP $(F,\alpha,\mu)$ is {\it regular} and we will denote by
$\omega_\mu$ the induced symplectic form on $F_\mu$. We agree to call $F_\mu$ the
symplectic quotient of  $(F,\alpha,\mu)$ even when this SFP  is not regular. 

The concept "symplectic factorization problem" is the symplectic analogon of the
concept "linearized action of a reductive group" on a polarized algebraic variety
in classical Geometric Invariant Theory; the relation between the two concepts is
beautifully explained in  [Kir].

 A {\it
compatible almost complex structure} of a symplectic factorization problem
$(F,\alpha,\mu)$ is an almost complex structure $J$ on $F$ which is
$\omega$-{\it tame} (which means that $\omega(\cdot, J\cdot)$ is a Riemannian
metric on $F$) and $K$-invariant.  Such an almost complex structure defines   an
$\omega_\mu$-tame complex structure $J_\mu$ on $F_\mu$, if the chosen SFP
was regular.

Many remarkable symplectic manifolds (e. g. Grassmann manifolds, flag manifolds,
toric varieties, etc) can be regarded in a natural way as symplectic quotients associated
with certain SFP's.  In many cases, the input symplectic manifold $F$ has -- in a
natural way -- a larger symmetry than the symmetry used in performing the
symplectic factorization. This larger symmetry induces then a symmetry on the
resulting symplectic quotient $F_\mu$, and plays an important role  in studying the
geometry of this quotient.   

\begin{dt} A symplectic factorization problem  with additional symmetry  (SFPAS)
is a 4-tuple
$(F,\alpha,K,\mu)$, where:
\begin{enumerate}
\item $F$ is a symplectic manifold, 
\item $\alpha:\hat K\times F\ra F$ is  an action of a compact Lie group
$\hat K$ on $F$, 
\item $K$ is a closed normal subgroup of $\hat K$, which acts symplectically on $F$ via
$\alpha$,
\item $\mu$ is a $\hat K$-equivariant moment map for the $K$-action on $F$.
\end{enumerate} 

A compatible almost complex structure of a SFPAS $(F,\alpha,K,\mu)$ is an almost complex
structure on $F$ which is $\omega$-tame and $\hat K$-invariant.
\end{dt}

In all interesting examples we know, $\alpha:\hat K\times F\ra F$ is itself
symplectic, and
$\mu$ is induced by a moment map for this $\hat K$-action via the projection
$\hat\kg^\vee\ra\kg^\vee$. The importance of these concepts comes from the
following obvious\\
\\
{\bf Remark:} {\it If $(F,\alpha,K,\mu)$ is a SFPAS, then the $\hat K$-action of $F$
descends to a   $K_0:=\hat K/K$-action on the symplectic quotient
$F_\mu$.  This action is symplectic if $\alpha$ was symplectic.

If $J$ is a compatible almost complex structure of the SFPAS
$(F,\alpha,K,\mu)$  and the SFP $(F,\alpha|_{K\times F},\mu)$ is regular, then the
induced $\omega_\mu$-tame almost complex structure $J_\mu$ on $F_\mu$  will
be $K_0$-invariant.}\\

Therefore a  SFPAS  (endowed with a compatible almost complex structure) provides
a symplectic quotient  (respectively an almost K\"ahlerian symplectic quotient)
which   comes with natural induced symmetry.

Below we give several relevant   examples, which motivate  the introduction of
these concepts and demonstrate their importance: 

\paragraph {1. Grassmann manifolds} \hfill{\break}

Consider the manifold $F=\Hom(\C^r,\C^{r_0})$ and the exact sequence of
compact Lie groups
$$1\map \U(r)\map \U(r)\times \U(r_0)\map \U(r_0)\map 1
$$
The  group $\hat K:=\U(r)\times\U(r_0)$ acts on $F$ in the obvious way.  We
denote by
$\alpha_{\rm can}$ this action. The
$K:=\U(r)$-action on $F$ has a one parameter family $(\mu_t)_{t\in\R}$ of
moment maps
$$\mu_t(f)=\frac{i}{2}f^*\circ f-it\id_{\C^r}\ ,
$$
and the corresponding symplectic quotients are
$$F_{\mu_t}=\left\{
\begin{array}{ccc}
\G r_{r}(\C^{r_0})&{\rm if}& t>0\\
\{*\}&{\rm if}& t=0\\
\emptyset&{\rm if}& t<0\ .
\end{array}\right.
$$

The moment maps $\mu_t$ are $\hat K$-equivariant, therefore the 4-tuples
$$(\Hom(\C^r,\C^{r_0}),\alpha_{\rm can},\U(r),\mu_t)$$
 are SFPAS's. This  implies    
that the $\hat K$-action on $F$ descends to a $K_0:=\U(r_0)$-action on the
symplectic quotients $F_{\mu_t}$.  This gives precisely the obvious
$\U(r_0)$-symmetry of the Grassamnnian $\G r_{r}(\C^{r_0})$ of $r$-planes of
the complex vector space $\C^{r_0}$.  This induced symmetry is essential for
understanding the geometry of the Grassmann manifolds.

\paragraph{2. Flag manifolds} \hfill{\break}
 
Let $V_1,\dots,V_m, V=V_{m+1}$ be Hermitian vector spaces. Put 
$$d_i:=\dim(V_i)\ ,\ d:=\dim(V)\ ,\ F:=\bigoplus_{i=1}^m \Hom(V_i,V_{i+1})\ .
$$
We consider the exact sequence of compact Lie groups
$$1\map\prod_{i=1}^m \U(V_i)\map \prod_{i=1}^{m+1} \U(V_i)\map \U(V)\map
1\ .
$$
and we put $K:=\prod_{i=1}^m \U(V_i)$, $\hat K:=\prod_{i=1}^{m+1} U(V_i)$,
$K_0:=\U(V)$.  The group $\hat K$ acts on $F$ by
$$\alpha_{\rm can}(g_1,\dots,g_{m+1})(f_1, \dots,f_m)=   (g_2\circ f_1\circ
g_1^{-1},  \dots,g_{m+1}\circ f_m\circ g_{m}^{-1})\ .
$$
The general form of a moment map for the restricted  $K$-action on $F$ is
$$\mu_t(f_1,\dots,f_m)=\frac{i}{2}
\left(
\matrix{f_1^*\circ f_1 \cr 
f_2^*\circ f_2-f_1\circ f_1^*\cr 
\dots\cr  
f_m^*\circ f_{m}-f_{m-1}\circ f_{m-1}^* }
\right)
-i\left(\matrix{t_1\id_{V_1}\cr
t_2\id_{V_2}\cr\dots\cr t_m\id_{V_m}}\right)
$$
where $t\in\R^m$.  To every $f=(f_1,\dots f_m)\in F$ we associate the
subspaces 
$$W_i(f):=(f_m\circ\dots\circ f_i)(V_i)\subset V\ ,\ 1\leq i\leq m$$
One obviously has $W_i\subset W_{i+1}$ and the map
$$f\mapsto (W_i(f))_{1\leq i\leq m}
$$
is constant on  orbits.
\begin{pr} Suppose that   $t_i>0$, for all $1\leq i\leq m$.
\begin{enumerate} 
\item Let $f\in F$.
Then the following conditions are equivalent:
\begin{enumerate}
\item $f$ is $\mu_t$-semistable
\item $f$ is $\mu_t$-stable
\item all maps $f_i$ are injective.
\end{enumerate}
\item The map
$$w:f\mapsto (W_i(f))_{1\leq i\leq m}
$$
identifies the symplectic quotient $F_{\mu_t}$ with the flag manifold
$$\F_{d_1,\dots d_m}(V):=\{(W_1,\dots,W_m)|\ W_1\subset
 \dots\subset W_m\subset V,\ \dim(W_i)=d_i\}$$
\end{enumerate}
\end{pr}

Note that   $\F_{d_1,\dots d_m}(V)$ is non-empty if and only if  $d_1\leq
 \dots\leq d_m\leq d$ and is interesting when all inequalities are strict (otherwise
it can be identified with a flag manifold associated with a smaller number of
$d_i$-s). We include a short proof for completeness.
\\ \\
\pf  1. We use the standard analytic criterion [B], [MU1], [LT] for testing stability.  

  If, in general, $a$ is a Hermitian endomorphism of a Hermitian
space $W$, and $\lambda\in \R$, we put
$$W^\lambda:=\ker(a-\lambda\id_W)\ ,\ W_\lambda:=\bigoplus_{\lambda'\leq
\lambda} W^{\lambda'}\ .
$$
 Let
$$\xi=(\xi_1,\dots,\xi_m)\in i\kg=\bigoplus_{i=1}^m \Herm(V_i)\ ,$$
and let $\alpha_*(\xi)\in\Herm(F)$ be the induced Hermitian endomorphism.  The
eigenvalue decomposition of this endomorphism is
$$\alpha_*(\xi)=\bigoplus_{i=1}^{m-1}
\sum_{\matrix{\scriptsize\eta\in Spec(\xi_{i+1}) 
\cr\scriptsize\lambda\in Spec(\xi_i)}}(\eta-\lambda)\id_{\Hom(V_i^\lambda,
V_{i+1}^\eta)}\oplus \sum_{\lambda\in
{\rm Spec}(\xi_m)}-\lambda\id_{\Hom(V_m^\lambda,V)}\ .
$$
A vector $f\in F$ is $\mu_t$-(semi)stable if and only if for every $\xi\in
i\kg\setminus\{0\}$ for which
$$f\in \bigoplus_{y\leq 0}{\rm Eig}(\alpha_*(\xi),y)\ ,
$$
one has 
$$\langle t,\xi\rangle=\sum_{i=1}^m   t_i\tr(\xi_i)> 0\ (\geq 0)\ .
$$
The condition $f\in \bigoplus_{y\leq 0}{\rm Eig}(\alpha_*(\xi),y)$ becomes
$$f_i(V_{i,\lambda})\subset V_{i+1,\lambda}\ \   \forall i\in \{1,\dots, m-1\},  \
\forall
\lambda\in{\rm Spec}(\xi_i)\ ,
$$
$$ f_m(V_{m,0})=\{0\}\ .
$$
If all $f_i$ are injective,   then this condition implies
$V_{m,0}=\{0\}$, and by induction, $V_{i,0}=\{0\}$  for all $i\in \{1,\dots,
m-1\}$. This means that all $\xi_i$ have only strictly positive eigen-values, hence
indeed $\langle t,\xi\rangle>0$.

Conversely, suppose that $f_i$ was not injective. Choose $\xi$ such that 
$\xi_j=0$ for all $j\ne i$,  and 
$$\xi_i=-{\rm pr}_{\id_{\ker (f_i)}}
$$
Then the condition  $f\in \bigoplus_{y\leq 0}{\rm Eig}(\alpha_*(\xi),y)$ is
obviously satisfied, but 
$$\langle t, \xi\rangle=-t_i\dim(\ker (f_i))< 0\ ,$$
hence $f$ is not $\mu_t$-semistable.\\

2.   The standard theory of K\"ahler quotients [Kir] gives
$$F_{\mu_t}=\qmod{\mu_t^{-1}(0)}{K}=\qmod{F_{\mu_t}^{st}
}{K^\C}=\qmod{\{f\in F| \ker(f_i)=0\ \forall i\in \{1,\dots,
m\}}{\prod_{i=1}^mGL(V_i)}\ ,
$$
where $F_{\mu_t}^{st}$ stands for the set of $\mu_t$-stable points in $f$.  It
easy to see that $\prod_{i=1}^mGL(V_i)$ acts freely on $F_{\mu_t}^{st}$, and
that the map $w$ identifies the quotient with $\F_{d_1,\dots d_m}(V)$.

\qed

Again, our input space $F$ has a larger symmetry that the $K$-symmetry used to
perform the symplectic factorization, namely the
symmetry defined by the 
$\hat K$-action.  The moment maps $\mu_t$ are all $\hat K$-equivariant, so  the 4-tuples 
$$(\bigoplus_{i=1}^m \Hom(V_i,V_{i+1}),  \alpha_{\rm can}, \prod_{i=1}^m
\U(V_i),\mu_t)
$$
are SFPAS's.

We get  an induced $K_0$-action on the symplectic quotients $F_{\mu_t}$. In the
particular case of Proposition 1.1 (when $t_i>0$), this is just the tautological
$\U(V)$-action on
$\F_{d_1,\dots d_m}(V)$.

\paragraph{3. Toric varieties}\hfill{\break}

In this case we take $F=\C^r$ and consider an exact sequence of the form
$$1\map K_w\map [S^1]^r\textmap{w} [S^1]^m \map1\ ,
$$
where $w:[S^1]^r\map[S^1]^m$ is an epimorphism and $K_w:=\ker(w)$.  Let
$v\in\Hom(\Z^r,\Z^m)$ be  the morphism defined by the differential $d_e(w)$.  The group
$\hat K:=[S^1]^r$ acts on
$F$ in a natural way.

The general form of a moment map for the induced $K_w$-action on $\C^r$ is
$$\mu_\tg(z^1,\dots,z^r)=-\frac{i}{2}p_v(|z^1|^2,\dots,|z^r|^2)+i\tg\ ,\
\tg\in\coker((v\otimes\id_\R)^*)\ ,
$$
where $(v\otimes\id_\R)^*:\R_m\map\R_r$ is the adjoint of the linear map induced
by the morphism $v$, and $p_v$ is the canonical projection
$$p_v:\R_r\map
\coker((v\otimes\id_\R)^*)=\qmod{\R_r}{(v\otimes\id_\R)^*(\R_m)}\ .
$$
Here we have used the natural identification
$${\rm Lie}(K_w)^\vee=i\coker[(v\otimes\id_\R)^*] \ .
$$

The quotient $F_{\mu_\tg}$ is non-empty if and only if 
$$\tg\in p_v(\R_r^{\geq 0})\ ,\ \R_r^{\geq 0}:=\{(t_1,\dots,t_r)\in\R_r|\
t_i\geq 0\}\ .
$$

Suppose that the following two conditions are satisfied:\\ \\
{\bf P$_1$: } {\it For every $j\in\{1,\dots,r\}$, the column  $v_j\in
\Z^m$ of the integer matrix $v$ is primitive, i. e. it is a generator of the semigroup
$\Z^m\cap
\R_{\geq 0} v_j$.}
\\ \\
{\bf P$_2$: } {\it  $\R_r^{\geq 0}\cap \im(v^*)=\{0\}$ in  the dual 
space  
$\R_r$ of $\R^r$. }
\vspace{2mm}\\

If these conditions are satisfied, then $F_{\mu_\tg}$ is a projective toric variety
with (at most) orbifold singularities, for every $\tg\in p_v(\R_r^{\geq 0})$ which
is a regular value of the map $z\mapsto p_v(|z_1|^2,\dots,|z_r|^2)$.   In order
to explain this statement more clearly, let us  recall  (see [Bat], [Gi]) some simple definitions
and results in the theory of toric varieties:

Fix a subset $J\subset\{1,\dots,r\}$, and let $\Sigma$ be a {\it
complete,  simplicial} fan   in
$\R^m$ whose 1-skeleton  $\Sigma(1)$ is the set of
rays $\Sigma(1)=\{\R_{\geq 0} v_j |\ j\in J\}$.
 For any $a=(a_1,\dots,a_r)\in \R_r$ and any strictly convex polyhedral cone
$\sigma\in\Sigma$ we define the functional  $f_\sigma^a$ on  the linear span
$\langle\sigma\rangle$ by  
$$\langle f_\sigma^a, v_j\rangle=-a_j\ \ {\rm if}\ \R_{\geq 0} v_j\ {\rm
is\ a \ ray\ of\ }\sigma\ .
$$

We introduce the following convex subsets of $\coker[(v\otimes\id_\R)^*]$ (see
[OT3] for details):
$$K(\Sigma):=\{p_v(a)\ \mid a_i\geq 0, \ \  \langle f^a_\sigma,
v_j\rangle\geq -a_j\ \forall \sigma\in\Sigma\ ,\  \forall j\in
\{1,\dots,r\}\}
$$
$$K_0(\Sigma):=\{p_v(a)\in  K(\Sigma)\ \mid \langle f^a_\sigma,
v_j\rangle> -a_j\ \forall \sigma\in\Sigma\ ,\  \forall j\in
\{1,\dots,r\} $$
$$\ {\rm for\ which}\ \R_{\geq	0} v_j\ {\rm is\ not\ a\ face\
of\ }\sigma\}
$$

A classical result in the theory of toric varieties states that  every complete simplicial
fan
$\Sigma$ in
$\R^m$ with $\Sigma(1)\subset\{\R_{\geq	0}v_1,\dots,\R_{\geq	0}v_r\}$ 
defines an associated     compact
toric variety $X_\Sigma$ in the follwing way: consider first the open set of $\C^r$
$$U(\Sigma)=\{z\in\C^r|\  \exists\sigma\in\Sigma \ {\rm such\ that}\
z^j\ne0\ 
\forall   j\in\{1,\dots,r\}\ {\rm for\ which}  $$
$$\ \R_{\geq 0}v_j \ {\rm
is\ not\ a\ face\ of}\ \sigma  
\}\ .
$$
One   proves that there is a {\it geometric quotient}
$$X_\Sigma:=\catqmod{U(\Sigma)}{K_w^\C}
$$
 and   this quotient is a compact algebraic 
variety with a natural orbifold structure; it is projective  if and only if
$K_0(\Sigma)\ne\emptyset $.

\begin{thry} Let $\Sigma$ be a complete simplicial fan $\Sigma$ in $\R^m$
with
$$\Sigma(1)\subset\{\R_{\geq	0}v^1,\dots,\R_{\geq	0}v^r\}\ .$$ 
  For every $\tg\in K_0(\Sigma)$, the set of semistable points with
respect to the moment map $\mu_\tg$ coincides with the correponding
set of stable points, and the symplectic quotient
$\qmod{\mu^{-1}_\tg(0)}{K_w}$ can be identified as a complex orbifold
with the projective toric variety $X_\Sigma$.
\end{thry}

Note also that {\it any} regular value of the map $z\mapsto
p_v(|z_1|^2,\dots,|z_r|^2)$ belongs to a $K_0(\Sigma)$ for a suitable complete
simplicial fan $\Sigma$ in $\R^m$, so one gets a complete description of all
  sympletic quotients  of $\C^r$ by $K_w$ which correspond to regular values of
the standard moment map $\mu_0$.\\

Concluding, we note that the 4-tuples 
$$(\C^r,\alpha_{\rm can},K_w,\mu_\tg)
$$
are obviously SFPAS's.
According to our general principle, the $[S^1]^r$-action $\alpha_{\rm can}$ on $F$ induces
an $[S^1]^m$-action on every symplectic quotient $F_{\mu_\tg}$. The complexification of
this action  has  a dense orbit and
 plays a fundamental role in the study of toric varieties.

\paragraph{4. Str{\o}mme's triples and Quot spaces of trivial sheaves on $\P^1$}
\hfill{\break}

 Consider three
Hermitian vector spaces $U$, $V$, $W$   of dimensions   $u$,
$v:=u+r$ and $w$ respectively, where $u$,  $w$, $r$ are non-negative integers 
with $r\leq w$.

This time our input symplectic manifold is
$$F:=\Hom(U,V)^{\oplus 2}\oplus\Hom(W,V)\ .
$$
hence the space of diagrams of the form
$$
\begin{array}{ccccc}
 &\stackrel{k}{\longrightarrow}& &&\\
U& &V&\stackrel{m}{\longleftarrow}&W\ .\\
 &\stackrel{l}{\longrightarrow}& &&
\end{array}
$$
We consider the exact sequence of compact Lie groups
$$1\map \U(U)\times \U(V)\map \U(U)\times \U(V)\times\U(W)\map \U(W)\map 1
$$
and we let the group $\hat K:=\U(U)\times \U(V)\times\U(W)$ act on $F$ by
$$\alpha(a,b,c)(k,l,m)=(b\circ k\circ a^{-1}, b\circ l\circ
a^{-1},b\circ m\circ c^{-1})\ .
$$

The general form of a moment map for the induced action of $K:=\U(U)\times
\U(V)$ on $F$ is given by
$$\mu_{s,t}(k,l,m):=\frac{i}{2}(k^*\circ k+l^*\circ l ,-k\circ k^*-       
l\circ l^*-m\circ m^*)+i(-t\id_{U}, s\id_{V})
$$
for real parameters $s,\ t\in\R$.  We refer to [LOT] for the following results
\begin{thry}  A triple $(k,l,m)$ is $\mu_t$-(semi)stable if and only if  
for all subspaces $U_1\subset U$, $V_1\subset V$ one has:
\begin{enumerate}
\item If 
$$(U_1,V_1)\ne(0,0)\ ,\ k(U_1)+l(U_1)\subset V_1
$$
then $s\dim(V_1)> t\dim(U_1)$ ($s\dim(V_1)\geq t\dim(U_1)$).
\item If 
$$(U_1,V_1)\ne(U,V)\ ,\ k(U_1)+l(U_1)\subset V_1\ ,\
\im(m)\subset V_1
$$
then \ $t\dim(U/U_1)> s\dim(V/V_1)$ ($t\dim(U/U_1)\geq
s\dim(V/V_1)$).
\end{enumerate}
\end{thry}
\begin{thry} Suppose   $s=1+\varepsilon$, $t=1$,  with $\varepsilon>0$ sufficiently small. 
\begin{enumerate} 
\item The first (semi)stability condition is  
equivalent to the condition that  the matrix $xk+yl$ has  maximal rank $u$ for a general pair $(x,y)\in\C^2$.
\item The second (semi)stability condition is  
equivalent to the condition that  the matrix $[xk+yl|m]\in M_{v,u+w}$ has  
maximal rank $v$ for all $(x,y)\in\C^2\setminus\{(0,0)\}$.
\end{enumerate} 
\end{thry}

A triple $(k,l,m)$ satisfying the two non-degeneracy conditions of Theorem 1.4 will be called
{\it  a Str{\o}mme triple}. The set ${\cal S}(u,w,r)$ of Str{\o}mme triples is
obviously open. The   importance of these objects comes from the following
construction of Str{\o}mme, which provides a simple description of certain Quot
spaces on $\P^1$.  With every Str{\o}mme triple $(k,l,m)$ we associate the diagram
$$
\begin{array}{cccccccc}
&&& &{\cal O}_{\P^1}\otimes W_{\phantom{b_b} }\ \ \  \ \ \ \ \ \ \ \  &&&
\\
&&&&\downarrow\overline m   \ \ \ \searrow \matrix{q_{k,l,m}\cr \phantom{b}}
\\
0&\ra&{\cal O}_{\P^1}(-1)\otimes U&\textmap{\overline{(k,l)}} &{\cal O}_{\P^1}\otimes V\
\map\
\   Q\ &\map 0\ ,
\end{array}
$$
where $\overline{(k,l)}$ is the sheaf morphism induced by $(k,l)$ and the standard isomorphism 
$$H^0({\cal O}_{\P^1}(1))\simeq\C^2\ ,$$
$\overline m$ is the sheaf morphism induced by $m$, and $q_{k,l,m}$ is the morphism which
makes the diagram     commute.  The fact that $(k,l,m)$ is a Str{\o}mme triple
means that
$\overline{(k,l)}$ is a monomorphism and
$q_{k,l,m}$  is an epimorphism.

An important result of Str{\o}mme [S] states that
\begin{thry}  The group $GL(U)\times GL(V)$ operates freely on the set of
Str{\o}mme triples and the map    $(k,l,m)\mapsto q_{k,l,m}$ descends to an
isomorphism
$$\qmod{{\cal S}(u,w,r)}{GL(U)\times GL(V)}\simeq Quot_{{\cal
O}_{\P^1}\otimes W}^{(r,u)}\ ,
$$
where $Quot_{{\cal O}_{\P^1}\otimes W}^{(r,u)}$ denotes the Quot space  classifying the
quotients of ${\cal O}_{\P^1}\otimes W$ which have rank $r$
and degree $u$.
\end{thry}

Our theorems Theorem 1.3, Theorem 1.4 imply that ${\cal S}(u,w,r)$ is just the
open set of $\mu_{1+\varepsilon,t}$ -(semi)stable points of $F$. Therefore, using
the analytic stability criterion, we see that the Quot space $Quot_{{\cal
O}_{\P^1}\otimes W}^{(r,u)}$ can be further
identified with a symplectic quotient:
\begin{co} \ [LOT]
One has a natural isomorphism
$$F_{\mu_{1+\varepsilon,1}}\simeq Quot_{{\cal O}_{\P^1}\otimes W}^{(r,u)}\ .
$$
\end{co}

Note finally that, in this case too, the input space comes with a larger symmetry than the  
$K$-symmetry
used to perform the symplectic factorization, namely with a natural
 $\hat K:=\U(U)\times\U(V)\times \U(W)$-symmetry. The moment maps $\mu_{s,t}$
are all $\hat K$-equivariant, so we get the following new  examples of SFPAS's:
$$
(\Hom(U,V)^{\oplus 2}\oplus\Hom(W,V), \alpha_{\rm can},\U(U)\times\U(V),\mu_{s,t})\ .
$$

Therefore, the symplectic quotients $F_{\mu_{s,t}}$ have all a natural
induced $U(W)$-symmetry. In the special case
$(s,t)=(1+\varepsilon,1)$, one gets precisely the obvious $U(W)$-action on the Quot space
$Quot_{{\cal O}_{\P^1}\otimes W}^{(r,u)}$.

\paragraph{5.  Quiver Problems}\hfill{\break}

All the examples above are particular cases of the following   {\it quiver
factorization problem}.

Let $Q=(V,A,s,t)$ be a quiver.  This means that $V$ and $A$ are finite sets (the set of
vertices and the set of arrows) and $s$, $t$ are maps $s,t:A\ra V$ (the
source map and the target map).

Let $\hat K$ be a compact Lie group. A $Q$-representation of $\hat K$ is a system 
$\rho=(\rho_v)_{v\in V}$, $\rho_v:\hat K\map U(W_v)$ of unitary representations of   
$\hat K$.

Take
$$F:=\bigoplus_{a\in A}\Hom(W_{s(a)},W_{t(a)})\ ,
$$
and consider the  $\hat K$-action $\alpha_\rho$ of $\hat K$ on $F$ induced by the  
$Q$-representation   $\rho$.  Let  
$K$ be a closed normal subgroup of
$\hat K$ endowed with an $\ad_{\hat K}$-invariant inner product on its Lie algebra $\kg$. 
The general form of a moment map for the induced
$K$-action is
$$\mu_\tg((f_a)_{a\in A})=\frac{i}{2}\rho^*_{i\kg}\left( (\sum_{s(a)=v} f_a^*\circ f_a
 -\sum_{t(a)=v}f_a\circ f_a^*)_{a\in A}\right)- i \tg\ ,
$$
where $\rho_{i\kg}$ is the composition 
$$\i\kg\hookrightarrow i\hat\kg\ \textmap{(id_e(\rho_v))_{v\in V}
}\ \bigoplus_{v\in V} iu(W_v)\ ,
$$
$\rho^*_{i\kg}$ is its adjoint, and $i\tg$ is a central element of $\kg$.

The 4-tuple 
$$(\bigoplus_{a\in A}\Hom(W_{s(a)},W_{t(a)}),\alpha_\rho,K,\mu_\tg)
$$
is obviously a SFPAS, so one gets an induced $K_0:=\hat K/K$-symmetry on the symplectic
quotients $F_{\mu_\tg}$. Such a SFPAS will be called a quiver factorization problem
associated with $Q$.\\
\\
{\bf Example:} Usually, one takes 
$$\hat K=\prod_{v\in V} U(W_v)\ ,$$
(where $W_v$ is Hermitian vector space) with the canonical representation  on $W_v$,  
and $K$ is chosen to be a closed normal subgroup of this product, for instance a  product
of factors $U(W_v)$.

In the very special case when $K =\hat K=\prod_{v\in V} U(W_v)$, one has a {\it
standard quiver factorization problem}. This is the case considered in classical GIT
[K].

The list below shows that all examples above are just special   quiver
factorization problems: \\ \\ \\
\hspace*{-30mm}
$
\begin{array}{|c|c|c|c|}
\hline
SFPAS&Q&\hat K& K\\
\hline
&&&\\
(\Hom(\C^r,\C^{r_0}),\alpha_{\rm
can},\U(r),\mu_t)&\bullet\map\bullet&\U(r)\times\U(r_0)&\U(r)\\
&&& \\
\hline
&&&\\
(\bigoplus\limits_{i=1}^m \Hom(V_i,V_{i+1}),  \alpha_{\rm can}, \prod\limits_{i=1}^m
\U(V_i),\mu_t)&
\matrix{\underbrace{\bullet\ra\bullet\ra\dots\ra\bullet}\cr
m+1}
&\prod\limits_{i=1}^{m+1}
\U(V_i)&\prod\limits_{i=1}^m
\U(V_i)\\
&&& \\
\hline
&&&\\
(\Hom(\C,\C)^{\oplus r},\alpha_{\rm can},K_w,\mu_\tg)&
\left.
\unitlength=1mm
\begin{picture}(9,12)(1,0)
\put(-6,0){$\bullet$}
\put(-3,3){\vector(1,1){6}}
\put(5,10){$\bullet$}
\put(-3,2){\vector(2,1){6}}
\put(5,5){$\bullet$}
\multiput(5.5,3)(0,-2.3){3}{.}
\put(-3,0){\vector(2,-1){6}}
\put(5,-5){$\bullet$}
\put(-3,-1){\vector(1,-1){6}}
\put(5,-10){$\bullet$}
\put(15,0){$r$}
\end{picture}
\right\} &[S^1]^r&K_w\\
&&& \\
\hline
&&&\\
(\Hom(U,V)^{\oplus 2}\oplus\Hom(W,V), \alpha,\U(U)\times\U(V),\mu_{s,t})
&\bullet\ \matrix{\ra\cr\ra}\ \bullet\ \leftarrow\ \bullet
&\U(U)\times\U(V)\times U(W)&\U(U)\times\U(V)
\\
&&& \\
\hline
\end{array}
$ 
\\
\\ \\ \\
{\bf Remark:} One can formulate a more general version of the quiver  factorization
problem explained above in the following way:

   Let
$\rho=(\rho_v)_{v\in V}$ be a
$Q$-representation of
$\hat K$, and let
$$r=(r_a)_{a\in A}\ ,\ r_a:\hat K\ra U(W_a^0)$$
 be a  system of representations of $\hat K$
indexed by the arrows of $Q$ (the "twisting representations").  Put 
$$F:= \bigoplus_{a\in A}\Hom(W_{s(a)},W_{t(a)}\otimes W_a^0)Ê,
$$
endowed with the $\hat K$-action $\alpha_{\rho,r}$ induced by $\rho$ and $r$, and let $K$
be a closed normal subgroup of $\hat K$.  The general form    of the moment map for the
induced
$K$-action on $F$ is 
$$\mu_\tg((f_a)_{a\in A})=\frac{i}{2}\rho^*_{i\kg}\left( (\sum_{s(a)=v} f_a^*\circ f_a
 -\sum_{t(a)=v}\tr_{W_v^0}(f_a\circ f_a^*))_{a\in A}\right)- i \tg\ , \ \i\tg\in z(\kg)\ .
$$
The system $(F,\alpha_{\rho,r},K,\mu_\tg)$ is obviously a SFPAS.  Such a SFPAS will be
called {\it an $r$-twisted  quiver factorization problem associated with $Q$}. The particular
case  when 
$$\hat K:=\prod_{v\in V} \U(W_v)\times \prod_{a\in A}\U(W_a^0)\ ,\   K:= 
\prod_{v\in V} \U(W_v)
$$
and $\rho_v$, $r_a$ are the canonical representations of $\hat K$ in $W_v$, $W_a^0$ was
considered in [AlPr1], [AlPr2].

\section{ The gauge theoretical problem associated with a symplectic factorization
problem with additional symmetry}

Let $F=(F,\omega)$ be symplectic manifold and let $\sigma_\mu=(F,\alpha,
,K,\mu)$ be a SFPAS.  We fix the following {\it  topological data}:
\begin{itemize} 
\item{ a closed, connected, oriented   real surface $Y$},
\item a principal $\hat K$-bundle $\hat P$ on $Y$,
\item a homotopy class $H\subset \Gamma(Y,E)$ of sections  in the associated
bundle
$$E:=\hat P\times_{\hat K} F\ .
$$
\end{itemize}

 To formulate our gauge theoretical problem, we also
need three {\it continuous parameters}:
\begin{itemize}
\item a compatible almost complex structure $J$  of the
SFPAS  $\sigma_\mu$, i. e.   a $\hat K$-invariant, $\omega$-tame almost
complex structure $J$ on the symplectic manifold $(F,\omega)$.
\item a Riemannian metric $g$ on $Y$,
\item  a parameter connection $A_0\in {\cal A}(P_0)$, where $P_0:=\hat P/K$ is
the associated $K_0:=\hat K/K$-bundle of $\hat P$.

\end{itemize}

We will denote by $J_g$ the complex structure defined by the Riemannian
metric $g$ and the fixed orientation of $Y$.

Note that any connection $\hat A$ on $\hat P$ defines an almost complex
structure $J_{\hat A}$ on $E$; this is the unique almost complex structure which
agrees with $J$ on the vertical tangent spaces of $E$ and with the complex
structure $J_g$ on the $\hat A$-horizontal spaces.\\

With $\sigma_\mu$ and our set of data we associate the following objects:\\
\begin{itemize} 
\item a {\it configuration space}
$${\cal A}:={\cal A}_{A_0}(\hat P)\times H\ ,
$$
where ${\cal A}_{A_0}(\hat P)$ stands for the affine space of those   
connections on $\hat P$ which induce $A_0$ on $P_0$.
\item a {\it differential equation} of vortex type for the elements $(\hat A,\varphi)\in 
{\cal A}$:
$$
\left\{\begin{array}{ccc}
\varphi&{\rm is}&J_{\hat A}-{\rm holomorphic}\\
{\rm pr}_{\kg} \Lambda_g F_{\hat A}+\mu(\varphi)&=&0\ . 
\end{array}
\right. \eqno{(V)}
$$
\item a {\it gauge group}
$${\cal G}:=\Gamma(Y,{\rm Aut}_{P_0}(\hat P))=
\Gamma(Y,\hat P\times_{\Ad} K)\ , 
$$
acting on the configuration space ${\cal A}$ and   leaving the set of solutions of the  equation
$(V)$ invariant.
\end{itemize}

Therefore, once the SFPAS $\sigma_\mu$  is fixed, the equation $(V)$ depends on 
two  systems of data: 
\begin{itemize}
\item the topological data $\tau=(Y,\hat P,H)$,
\item the   continuous parameters $\pg:=(J,g,A_0)$\ .
\end{itemize}

When we want to take into account this dependence, we will write
$(V^\tau_\pg)$, or -- when $\tau$ is obvious -- $(V_\pg)$ instead of $(V)$.

We denote by ${\cal A}^*$ the open subset of ${\cal A}$ consisting of
{\it irreducible pairs}, i. e.   pairs
$(\hat A,\varphi)$ with trivial stabilizers with respect to the ${\cal G}$-action,
and by ${\cal A}^V$ ($[{\cal A}^*]^V$) the closed subspace of ${\cal A}$
(${\cal A}^*$) of solutions of $(V)$.  

We introduce the quotients
$${\cal B}:=\qmod{{\cal A}}{{\cal G}}\ ,\ {\cal B}^*:=\qmod{{\cal A}^*}{{\cal G}}
$$
and the
moduli spaces of solutions
$${\cal M}:=\qmod{{\cal A}^V}{{\cal G}}\subset {\cal B}\ ,\  {\cal
M}^*:=\qmod{[{\cal A}^*]^V}{{\cal G}}\subset {\cal B}^*\ .
$$
After suitable Sobolev completions, ${\cal B}^*$ becomes a Banach manifold,    
and, by standard  elliptic theory,  ${\cal M}^*$ becomes a finite dimensional  subspace of
this manifold. We will also use the notations ${\cal M}_\pg$,
${\cal M}^\tau_\pg$,
${\cal M}^\tau_\pg(\sigma_\mu)$ when we have to take into account the
functoriality of these objects. 

As in Donaldson theory or in classical Gromov-Witten theory, in order to introduce
invariants associated with these moduli spaces, we have to endow them with
{\it canonical cohomology classes} and with a {\it virtual fundamental class}.\\

Canonical cohomology classes on ${\cal M}^*$ can be obtained in the following
way:

We regard a section $\varphi\in \Gamma(Y,E)$ as a $\hat K$-equivariant map
$\hat P\ra F$.  We have a natural evaluation map
$${\rm ev}:{\cal A}\times \hat P\map F\ ,
$$
which is $\hat K$-equivariant and ${\cal G}$-invariant, hence  descends to a
$\hat K$-equivariant map 
$$\qmod{{\cal A}\times \hat P }{{\cal G}}\map F
$$
which restricts to a $\hat K$-equivariant map 
$$\Phi:\hat {\cal P}:=\qmod{{\cal A}^*\times \hat P }{{\cal G}}\map F\ .
$$
The space $\hat {\cal P}$ can be regarded as $\hat K$-bundle  over the product ${\cal
B}^*\times Y$ (the {\it universal $\hat K$-bundle}), whereas $\Phi$ can be
interpreted as {\it the universal section} in the associated bundle $\hat{\cal
P}\times_{\hat K}F$.  The map $\Phi$ induces a morphism
$$\Phi^*:H^*_{\hat K}(F)\map H^*({\cal B}^*\times Y)\ .
$$
For any cohomology class $c\in H^*_{\hat K}(F,\Z)$ and homology class $h\in
H_*(Y,\Z)$   put
$$\delta^c(h):=\Phi^*(c)/h\in H^*({\cal B}^*,\Z)\ .
$$

We will say that the pair $(\pg,\mu)$ is {\it good} if ${\cal
M}^\tau_\pg(\sigma_\mu)= {\cal M}^\tau_\pg(\sigma_\mu)^*$, i. e. the
corresponding equation $(V)$ has only irreducible solutions.

In many interesting cases (see  Theorem 2.2 below) one can show that ${\cal M}^*$ 
can be identified with the vanishing locus of a Fredholm section $\vg$ (induced by
the left hand term $v$ of $(V)$) in a Banach vector bundle over the Banach
manifold
${\cal B}^*$, and that the determinant line bundle of the index of the family of
intrinsic differentials of  $\vg$ can be naturally oriented in a neighbourhood of this
vanishing locus. In this case, the formalism of Brussee ([Bru], [OT2]) applies and yields
a virtual fundamental class  
$$[{\cal M}^*]^{\rm vir}\in H_{\rm index(\vg)}^{\rm cl}({\cal M}^*,\Z)\ .$$ 
in the homology   with closed supports of ${\cal M}^*$.

 If this situation occurs, and if the moduli space ${\cal M}$ is compact (see Theorem 2.2 below)
one can define {\it   gauge theoretical Gromov-Witten   invariants of the SFPAS
$\sigma_\mu$} with respect to the parameters
$(\tau,\pg)$ by
$$GGW^\tau_\pg(\sigma_\mu)\left(\left(
\matrix{c_1\cr h_1}\right),\dots,\left(\matrix{c_k\cr
h_k}\right)\right):=\left\langle \cup_{i=1}^c \delta ^{c_i}(h_i), [{\cal
M}^\tau_\pg(\sigma_\mu)]^{\rm vir}\right\rangle \ ,
 $$
for every good pair $(\pg,\mu)$.

\begin{re} The numbers $GGW^\tau_\pg(\sigma_\mu)\left(\left(
\matrix{c_1\cr h_1}\right),\dots,\left(\matrix{c_k\cr
h_k}\right)\right)$ are not independent, but they satsfy a set of \ub{tautological}
\ub{relations}. Therefore   the map $GGW^\tau_\pg(\sigma_\mu)$ descends to a
graded $\Z$-algebra $\A$ generated by the symbols $\left[ \matrix{c\cr h}\right]$,
$c\in H^*_{\hat K}(F;\Z)$, $h\in H_*(Y,\Z)$ subject to these tautological  relations 
[OT2]. \\

For instance, if $\hat K=\U(r)\times K_0$ and $F$ is contractible,  one has
$$\A=\Z[u_1,\dots,u_r,v_2,\dots,v_r]\otimes \Lambda^*[\oplus_{l=1}^r  H_1(Y,\Z)_l]\ ,
$$
$$\deg(u_i)=2i,\ \deg(v_j)= 2j-2,\ \deg(H_1(Y)_l)=2l-1\ ,
$$
$$u_i=\left(\matrix{{\rm chern}_i\cr
[*]}\right)\ , \ v_j=\left(\matrix{{\rm chern}_j\cr
[Y]}\right)\ ,\ H_1(Y,\Z)_l=\left\{\left.\left(\matrix{{\rm chern}_l\cr
h}\right)\right|h\in H_1(Y,\Z)\right\}\ .
$$
Here ${\rm chern}_l\in H^{2l}_{\hat K}(F,\Z)$ is induced by the isomorphism
$H^{*}_{\hat K}(F,\Z)\simeq H^*(B\hat K,\Z)$ and the natural map $H^*(B\hat K,\Z)\ra
H^*(B\U(r),\Z)$.
\end{re}

\begin{thry}\hfill{\break}
 1.  If the action $\alpha:\hat K \times F\ra F$ is defined by a unitary
representation of $\hat K$, then the moduli space ${\cal M}^*$ can be regarded
as the vanishing locus of a holomorphic Fredholm section $\vg$ in a complex Banach
vector bundle over the Banach manifold ${\cal B}^*$.\\
2. If $\alpha|_{K\times F}$ is defined by a unitary
representation $\rho:K\ra U(F)$ of $K$ and the standard moment map
$$\mu_\rho(f)= (d_e\rho)^*(-\frac{i}{2} f\otimes f^*)
$$
satisfies the properness condition $\mu_\rho^{-1}(0)=\{0\}$, then the moduli
space ${\cal M}^\tau_\pg(\sigma_\mu)$ is always compact  for all parameters
$(\tau,\pg)$ and any  moment map $\mu$.
\end{thry}

Therefore, in the first case   Brussee's formalism [Bru]
applies and gives a virtual fundamental class $[{\cal M}^*]^{\rm vir}\in H_{\rm
index(\vg)}^{\rm cl}({\cal M}^*,\Z)$, whereas in the second situation all the moduli spaces
${\cal M}^\tau_\pg(\sigma_\mu)$ are compact.
Combining these results it follows that
\begin{co} If $\alpha$ is defined by a unitary
representation
$\hat
\rho$ and
$\mu_{\hat \rho|_K}$  satisfies the properness condition above, then the gauge
theoretical Gromov-Witten invariants $GGW^\tau_\pg(\sigma_\mu)$ are well
defined for every good pair $(\pg,\mu)$.
\end{co}
\begin{re} One should be able to define $\Q$-valued invariants for \ub{almost}
\ub{good} pairs $(\pg,\mu)$, i. e. for pairs $(\pg,\mu)$ for which the solutions of
the corresponding equation $(V)$ have  only finite stabilizers. In this case, the
moduli space ${\cal M}$ must be endowed with an orbifold structure. We refer to
[CMS] and [OT3] for details concerning this generalization.
\end{re}

Let us explain   two important special cases of our gauge theoretical problem  
obtained by chosing $K$ trivial or $K_0$ trivial:
\begin{enumerate}
\item       $K=\{1\}$, $F$ compact:

In this case the second equation in the system $(V)$ is identically satisfied, $\hat K=K_0$,
$\hat P=P_0$, and the associated bundle $E$ has a fixed  almost complex
structure $J_{A_0}$ induced by $A_0$. In this case
$${\cal M}=\Gamma_{J_{A_0}}(Y,E) 
$$
is the space of $J_{A_0}$-almost holomorphic sections in $E$.  Of course, ${\cal
M}$ is in general non-compact, because bubbling phenomena occur, exactly as in
classical Gromov-Witten theory. However, generalizing ideas of Ruan [R] one can
define a natural compactification $\bar  {\cal M}$ of such a moduli space  and, at
least for special fibres $F$, define Gromov-Witten type invariants.

These Gromov-Witten invariants should be called $K_0$-{\it twisted
Gromov-Witten invariants}, because they are obtained by replacing the moduli
spaces of almost holomorphic $F$-valued maps   in
Ruan's Gromov-Witten theory by moduli spaces of sections in $F$-bundles  
with symmetry group
$K_0$ (i. e. moduli spaces of $K_0$-twisted
$F$-valued maps).
\item  $K_0=\{1\}$, $F_\mu$ compact: 

In this case $\hat K=K$, so one gets a moduli problem which depends on a SFP
$(F,\alpha,\mu)$ rather than a SFPAS.  The obtained invariants should be called
gauge theoretical (or Hamiltonian) $K$-{\it equivariant Gromov-Witten invariants}
of $F$.  This case was extensively considered and studied by   Mundet i Riera [Mu2],
Gaio [Ga], Gaio-Salamon [GS] Cieliebak-Gaio-Salamon [CGS], 
Cieliebak-Gaio-Mundet--Salamon [CGMS].  They state a conjecture -- called the {\it
adiabatic limit conjecture} -- asserting that {\it the standard Gromov-Witten
invariants of the symplectic quotient
$F_\mu$ can be expressed in terms of gauge theoretical
Gromov-Witten invariants of the SFP $(F,\alpha,\mu)$}, provided this SFP is regular, 
the quotient $F_\mu$ is compact, and both types of invariants are well defined.

This conjecture is motivated  by the fact that the  second
equation  in $(V)$ tends to 
$$\mu(\varphi)=0$$
 when $g$ is replaced by $tg$ and
$t\ra\infty$.  Progress on this problem   was obtained in [GS]
for a special class of symplectic manifolds $F$.
\end{enumerate}

Taking into account these discussions, one should probably think of the invariants
associated with a SFPAS $(F,\alpha,K,\mu)$ as {\it  gauge theoretical $K_0$-twisted,
$K$-equivariant Gromov-Witten  invariants} of $F$.

Of course, the adiabatic limit conjecture should  also be true for the general  
twisted equivariant Gromov-Witten invariants, but one should replace  the usual
Gromov-Witten invariants of the symplectic quotient $F_\mu$ by its
$K_0$-twisted invariants.

\section{The complex geometric interpretation}

The natural question at this point is how can one obtain explicit descriptions
of moduli spaces ${\cal M}^\tau_\pg(\sigma_\mu)^*$ (${\cal
M}^\tau_\pg(\sigma_\mu)$) associated with a given SFPAS $\sigma_\mu$ and
parameters $\tau=(Y,\hat P,H)$,
$\pg=(J,g,A_0)$.\\

The main observation here is \\

{\it When the SFPAS $\sigma_\mu$ is \ub{K\"ahlerian}, i. e. the almost complex structure $J$
on
$F$    {is}  {integrable} and the
$\hat K$-action on $F$ extends to a holomorphic $\hat K^\C$-action, then the   
\ub{universal}
\ub{Kobayashi}-\ub{Hitchin} \ub{correspondence} gives a complex geometric
interpretation of these moduli spaces in terms of (poly)\ub{stable}  \ub{framed}  
 \ub{holomorphic} \ub{pairs}.
}\\

There are several concepts here which must be explained: \\

We denote by $G$, $\hat G$ and $G_0$ the complexifications of $K$, $\hat K$
and $K_0$; these are complex  reductive groups.  Let
$\hat Q$ ($Q_0$) be the $\hat G$-bundle ($G_0$-bundle) obtained by complexifying 
$\hat P$ ($P_0$).  The connection $A_0$ induces  a bundle-holomorphic
structure $J_0$ on $Q_0$\footnote{We recall that, in general,  a bundle-holomorphic
structure on a principal $G$-bundle $Q$ is a  holomorphic structure on the total space with
respect to which the action $G\times Q\ra Q$ of the structure group and the projection on the
base are both holomorphic.}, so one obtains a holomorphic $G_0$-bundle ${\cal Q}_0$ on the
Riemann surface $(Y,J_g)$.

\begin{dt} A {\it framed holomorphic pair  of type} $(Q,\alpha,H,{\cal Q}_0)$ is a
triple   
 $(\hat {\cal Q},\varphi,\lambda)$, where 
\begin{itemize}

\item $\hat {\cal Q}$ is  a holomorphic
$\hat G$-bundle over $Y$

\item $\varphi$ is a holomorphic section in the associated
$F$-bundle $\hat {\cal Q}\times_{\hat G} F=\hat P\times_{\hat K} F$ belonging
to $H$,
\item   $\lambda:\hat{\cal Q}\ra{\cal
Q}_0$ is $G_0$-{\it framing} of $\hat Q$, i. e. it is a holomorphic bundle morphism
of type $\hat G\ra G_0$,
\end{itemize}
such that $\hat {\cal Q}$ is ${\cal C}^\infty$-isomorphic to $\hat Q$ over
$Q_0$, i. e.  there exists a commutative diagram\\
$$
\begin{array}{c}
\unitlength=1mm
\begin{picture}(20,12)(-5,-6)
\put(-6,4){$\hat{\cal Q}$}
\put(0,5){\vector(2,0){8}}
\put(10,4){$\hat Q$}
\put(2,6){${\cal C}^\infty$}
\put(2,2){$\simeq$}
\put(-3,3){\vector(2, -3){4}}
\put(2,-6){${  Q}_0$}
\put(-4,-3){$\lambda$}
\put(10,3){\vector(-2, -3){4}}
\end{picture} 
\end{array}\ ,
$$
 where the right hand morphism is just the canonical projection $\hat Q\ra
Q_0$.
\end{dt}

An {\it isomorphism} between such triples $(\hat {\cal Q},\varphi,\lambda)$, $(\hat
{\cal Q}',\varphi',\lambda')$ is a bundle
isomorphism
$f:\hat {\cal Q}\ra \hat {\cal Q}'$ such that $\lambda'\circ f=\lambda$ and
$f_*(\varphi)=\varphi'$.\\

The classification of   framed holomorphic pairs of a given type is a very
interesting and   important complex geometric problem.  Many moduli problems
in complex geometry are special cases of this "universal" 
classication problem.  One can develop a complex geometric       {\it deformation
theory} for such objects, i. e. one can introduce in a natural way the notions of
{\it holomorphic families of framed holomorphic pairs} (of a fixed type)     
 parameterized by a complex space, {\it versal} and {\it universal} deformation of a fixed
framed pair, etc.

As in classical GIT, one cannot construct a moduli space with good properties
classifying all framed holomorphic  pairs of a given type; one needs a {\it stability
condition}.

The point is that with every pair $(\mu,g)$ consisting of a $\hat K$-equivariant 
moment map for the $K$-action on $F$ and a Riemannian metric on $Y$ one can
associate a stability    condition. The stability condition is open, 
i. e. if $(\hat {\cal Q}_t,\varphi_t,\lambda_t)_{t\in T}$ is a holomorphic family of
framed holomorphic pairs parameterized by a complex space $T$, then the set of
parameters $t$ for which  $(\hat {\cal Q}_t,\varphi_t,\lambda_t)$ is
$(\mu,g)$-stable is open in
$T$.  Moreover, using complex geometric deformation theory, one gets a   moduli
space ${\cal M}^{(\mu,g)-{\rm st
}}(Q,\alpha,H,{\cal Q}_0)$ classifying $(\mu,g)$-stable framed pairs of type
$(Q,\alpha,H,{\cal Q}_0)$.

Writing  the  stability conditions explicitely in the general case is quite technical and
requires a long preparation, but is now well understood (see [Mu1] for the case $K_0=\{1\}$,
[LT] for the general case).  There exist two important special cases in which it takes a simple
form:
\begin{pr} Suppose that $K$ is abelian and that one of the following conditions is
satisfied:
\begin{enumerate}
\item $F$ is a Hermitian space and $\alpha|_{K\times M}:K\times F\ra F$ is induced   
by a unitary representation, or
\item $F$ is a quasiprojective, and $(\omega,J)$ are   induced by a regular embedding 
$F\ra\P^N$,  and $\alpha|_{K\times M}:K\times F\ra F$ is induced    by a unitary
representation in $\C^{N+1}$.
\end{enumerate}

Then $(\hat{\cal Q},\varphi,\lambda)$ is $(\mu,g)$-stable if and only if $\varphi$ is
generically
$(\mu-\frac{2\pi}{Vol_g(Y)}\mu_K(\hat Q))$-stable.
\end{pr}

Here $\mu_K(\hat Q)$ denotes a topological invariant of $\hat Q$ which generalizes the usual
slope of a vector bundle [LT]. 

The universal Kobayashi-Hitchin correspondence states that, with the notations and
conventions above, there are natural isomorphisms
$$
\begin{array}{c}{\cal M}^\tau_\pg(\sigma_\mu)^*\simeq {\cal M}^{(\mu,g)-{\rm st
}}(Q,\alpha,H,{\cal Q}_0)  \ .
\end{array}\eqno(KH)
$$
In particular, this show that ${\cal M}^\tau_\pg(\sigma_\mu)^*$
has a natural complex space structure.
\begin{re} \hfill{\break} 
1. The Kobayashi-Hitchin correspondence can be extended to
arbitrary compact complex manifolds $Y$.  When $\dim_\C(Y)\geq 2$, one must add
the integrability condition $F^{02}_{\hat A}=0$ to the system $(V)$ and require that
the parameter connection $A_0$ is integrable.\\
2. There is a more refined Kobayashi-Hitchin correspondence which identifies the whole
moduli space ${\cal M}^\tau_\pg(\sigma_\mu)$ with the moduli space of
$(\mu,g)$-\ub{polystable} framed pairs [LT].
\end{re}

The Kobayashi-Hitchin correspondence is a very important tool which can be used to
give explicit descriptions of moduli spaces, but    it is not quite sufficient for the
computation of the invariants.  In order to have a complex geometric
interpretation of the invariants, one also has to compare {\it the virtual
fundamental classes} of the moduli spaces involved in the Kobayashi-Hitchin
correspondence $(KH)$.  As we explained in the previous section, a moduli space
${\cal M}^*$ defined with gauge theoretical methods can  naturally be endowed  
with a virtual fundamental class $[{\cal M}^*]^{\rm vir}$ when it can be identified
with the vanishing locus of a Fredholm section in a vector bundle over a Banach  
manifold, and the determinant line bundle of the index of the family of intrinsic
differentials of this section is oriented in a neighbourhood of this vanishing locus. 
We agree to call such a gauge theoretical moduli problems {\it of Fredholm type}.

On the other hand, many interesting moduli spaces defined in Algebraic Geometry
come with a natural {\it perfect obstruction theory} in  the sense of
Behrend-Fantechi [BF], and such a structure allows one to define a virtual fundamental
class in the Chow group $A_*$ of the moduli space ([BF],  [Kr]).

We state that\\ \\
{\bf Conjecture}: {\it Let 
$$\iota: {\cal M}^*\textmap{\simeq}{\cal M}^{\rm st}
$$
be any Kobayash-Hitchin type isomorphism. Suppose that

\begin{itemize}
\item ${\cal M}^*$ is associated with a gauge theoretical moduli problem of Fredhom
type,
\item all the data involved in the definition of ${\cal M}^{\rm st}$ are algebraic.
\end{itemize}

Then ${\cal M}^{\rm st}$ admits a canonical perfect obstruction theory, and  $\iota$ maps
the
 virtual fundamental class of ${\cal M}^*$ in the sense of Brussee to the image  of the
Behrend-Fantechi virtual fundamental class of ${\cal M}^{\rm st}$ 
 under the cycle map 
$$cl: A_{k}({\cal M}^{\rm st}) \ra H_{2k}^{\rm cl}({\cal M}^{\rm
st},\Z)\ .$$
}
\begin{re} \hfill{\break} 
1. The second assumption can  probably be removed, but in order to give a
sense to this more general statement, one   needs a complex analytic generalization of the
Behrend-Fantechi virtual class theory.\\
2. The moduli problem introduced in Chapter 2 is always of Fredholm type when the action
$\alpha$   is  defined by a unitary representation.  Note however that the analogous vortex
type problems on higher dimensional base manifolds $Y$ are in general \ub{not}  of Fredholm
type.\\ 3. The standard vortex moduli problem on K\"ahler surfaces can be regarded as a
moduli problem of Fredholm type [OT2].
\end{re}

The importance of our conjecture is obvious: it provides a universal principle allowing one to
identify not only certain moduli spaces defined within the two theories, but also the
corresponding numerical invariants. 

The conjecture has already been checked  in several special cases [OT2], [OT3], and we
believe that   our method  can be used to give a  proof of it under the
assumption that   the base manifold is a Riemann surface.

\section{Examples, computations, and applications}

We illustrate the universal Kobayashi-Hitchin correspondence in some of the special
situations  which we  considered in the first chapter of this article.  In every case we will     
indicate explicitely the corresponding stability condition and we will state several results 
concerning the corresponding gauge theoretical invariants, as well as  applications of
these results.

\paragraph{1. The SFPAS which yields the Grassman manifolds}\hfill{\break}

We come back to the SFPAS
$$\sigma_{\mu_t}=(\Hom(\C^r,\C^{r_0}),\alpha_{\rm can}, U(r),\mu_t)\ ,
$$
where   $\alpha_{\rm can}$ is the natural action of $\hat K:=U(r)\times U(r_0)$ on
$F=\Hom(\C^r,\C^{r_0})$  endowed with its obvious compatible complex 
structure $J$, and 
$$\mu_t(f)=\frac{i}{2}f^*\circ f-it\id_{\C^r}\ .
$$

The data of a principal $\hat K$-bundle $\hat P$ on a real surface $Y$ is equivalent
to the data of a pair
$(E,E_0)$ of Hermitian vector bundles  of ranks $r$, $r_0$ over $Y$.  Let $d$,
$d_0$ be the degrees of these bundles.  We denote by $\tau$ the topological data
$(Y,\hat P,H)$, where, of course, $H:=A^0\Hom(E,E_0)$.

Consider a Hermitian connection $A_0$ on $E_0$, and denote by
${\cal E}_0$ the corresponding holomorphic bundle.    We also fix a Riemannian
metric $g$ on $Y$.  So our continuous data is the triple $\pg=(J,A_0,g)$.

The corresponding   gauge theoretical problem becomes: For a given  real
number $t$, classify   pairs
$(A,\varphi)$ consisting of a  Hermitian connection in $E$ and a  morphism
$\varphi\in A^0\Hom(E,E_0)$ such that the following vortex type equation is
satisfied:
$$
\left\{\begin{array}{ccc}
\bar\partial_{A,A_0}\varphi&=&0\\
i\Lambda F_A-\frac{1}{2}\varphi^*\circ\varphi &=& -t\id_E\ .

\end{array}\right. \eqno{(V^{A_0}_t)}
$$
Let  ${\cal M}_t(E,E_0,A_0)$ be the moduli space of solutions of this system and 
${\cal M}^*_t(E,E_0,A_0)$   the open
subspace of irreducible solutions.\\

We explain now the   stability condition which corresponds to this gauge
theoretical problem  [Bra], [HL], [OT2]: 

Let $\tau$ be a real constant with $\deg(E)/\rk(E)>-\tau$.
A pair $({\cal E},\varphi)$ consisting of a holomorphic bundle ${\cal E}$ of ${\cal
C}^\infty$-type $E$ and a holomorphic   morphism ${\cal
E}\textmap{\varphi}{\cal E}_0$ is
$\tau$-(semi)stable  if  for every nontrival
   subsheaf ${\cal F}\subset {\cal E}$  one has
$$ \begin{array}{cccc} \mu({\cal E}/{\cal F})   &(\geq) &-\tau& {\rm if}\
\rk({\cal F})< r,\\
   \mu({\cal F})&(\leq)&-\tau  &Ö \ \ \ \ {\rm if }\ \ 
{\cal F}\subset\ker(\varphi)\ .
\end{array}
$$
Let ${\cal M}^{st}_\tau(E,{\cal E}_0)$ be the moduli space of $\tau$-stable holomorphic
pairs $({\cal E},\varphi)$ as above. Then the Kobayashi-Hitchin correspondence states:

\begin{thry} \ [Bra], [OT2]   There is a natural real analytic isomorphism 
$${\cal
M}^*_t(E,E_0,A_0)\textmap{\simeq}{\cal M}^{st}_\tau(E,{\cal E}_0)\ ,$$
where $\tau=\frac{  \Vol_g(Y)}{2\pi} t$.
\end{thry}

The complex geometric moduli spaces ${\cal M}^{st}_\tau(E,{\cal E}_0)$ are quite
complicated in general. However, if either $r=1$, or $\tau\gg 0$ they have a 
beautiful algebraic geometric interpretation in terms of Quot spaces.              In
general, for a ${\cal C}^\infty$-bundle $F$ and a holomorphic bundle ${\cal F}_0$
over $Y$ we   denote by $Quot^F_{{\cal F}_0}$ the Quot space of
quotients of ${\cal F}_0$ with kernels differentiably isomorphic to $F$.

\begin{pr} \hfill{\break} 
i)  Suppose $r=1$.  Then 
$${\cal M}_\tau^{\rm st}(E,{\cal E}_0)=\left\{\begin{array}{ccc}
\emptyset&{\rm if}& \tau <
- \frac{d}{r}\ \  \\
Quot^E_{{\cal E}_0}&{\rm if}& \tau >
- \frac{d}{r}\ .
\end{array}\right.
$$
ii) There exists a constant $c({\cal E}_0,E)$ such that  for all
$\tau\geq c({\cal E}_0,E)$ one has:
\begin{enumerate}
\item For every  $\tau$-semistable pair $({\cal E},\varphi)$,   $\varphi$
is injective. 
\item Every pair  $({\cal E},\varphi)$ with $\varphi$
injective   is $\tau$-stable. 
\item There is a natural isomorphism ${\cal M}^{st}_\tau(E,{\cal
E}_0)=Quot^E_{{\cal E}_0}$. 
\end{enumerate}
\end{pr}

Combining this result with Theorem 4.1, and taking into account the reducible solutions, one
obtains 
\begin{co}\hfill{\break}
i) In the abelian case   $r=1$   one has 
$${\cal
M}_t (E,E_0,A_0)={\cal
M}_t^*(E,E_0,A_0)\ \ {\rm for}\
 t\ne-\frac{2\pi}{\Vol_g(Y)}\frac{d}{r}\ ,$$  
and a real analytic isomorphism 
$${\cal M}_t(E,E_0,A_0)\simeq\left\{\begin{array}{ccc}
\emptyset&{\rm if}& t <
-\frac{2\pi}{\Vol_g(Y)}\frac{d}{r}\ \  \\
Quot^E_{{\cal E}_0}&{\rm if}& t >
-\frac{2\pi}{\Vol_g(Y)}\frac{d}{r}\ .
\end{array}\right.
$$
\\
ii) For sufficiently large $t\in\R$   one has  ${\cal
M}_t(E,E_0,A_0)={\cal M}^*_t(E,E_0,A_0)$ and a natural
identification
$${\cal M}_t(E,E_0,A_0)\simeq Quot^E_{{\cal E}_0} .$$
\end{co} 
\begin{re}  All these results can be generalized to the analoguous vortex problems on
arbitrary K\"ahler manifolds [OT2].
\end{re}

Using these results and well-known explicit descriptions of the 
Quot spaces on Riemann surfaces, we are able to compute explicitely all twisted
equivariant gauge theoretical Gromov-Witten invariants
$GGW^\tau_\pg(\sigma_{\mu_t})$ in the abelian case  $r=1$.  Using Remark 2.1,
put
$$GGW^\tau_\pg(\sigma_{\mu_t})(l):=\sum_{i=0}^\infty
  GGW^\tau_\pg(\sigma_{\mu_t})(\left(\matrix{{\rm chern}_1\cr
[*]}\right)^i\otimes \left(\matrix{{\rm chern}_1\cr
l}\right))\ ,
$$
for every $l\in \Lambda^*(H^1(Y,\Z))$.  Only one term in this series  	 can be
non-zero if $l$ is homogeneous.

\begin{thry} \ [OT2] Suppose $r=1$,   and let $v$ be the expected
dimension of the moduli space, i. e. 
$v:=d_0-r_0d+(r_0-1)(1-g(Y))$.
Let $l_{\ooo_1}$ be the generator of $\Lambda^{2g}(H^1(Y,\Z))$ defined
by the complex orientation $\oo_1$ of $H^1(Y,\R)\simeq H^{0,1}(Y)$. Then   
$$GGW_\pg^{\tau}(\sigma_{\mu_t})(l)=\left\{
\begin{array}{ccc}
\left\langle\sum\limits_{i\geq\max(0,g(Y)-v)}^{g(Y)}\frac{(r_0
\Theta)^i}{i!}\wedge l\ ,\ l_{\ooo_1}\right\rangle&{\rm if}& t >
-\frac{2\pi}{\Vol_g(Y)}\frac{d}{r}\\
0&{\rm if}& t <
-\frac{2\pi}{\Vol_g(Y)}\frac{d}{r}\ .
\end{array}
\right.
$$
\end{thry}
\vspace{2mm}

{\bf Applications:}\\ \\  
1. If one admits the adiabatic limit conjecture (see section 2) this result would allow
the full computation of all "{\it higher genus}" twisted   Gromov-Witten invariants of
$\P^n$.

Note that, according to some experts, the standard (Kontsevich) higher genus
Gromov-Witten invariants of $\P^n$ are not completely understood  (see for
instance the talk "Relative Gromov-Witten invariants" by  Andreas
Gathmann (Princeton), Oberwolfach 2002).\\
\\
2. Let $Quot^E_{{\cal E}_0}$ be a  Quot space of    dimension zero and
vanishing expected dimension $v=r d_0-r_0 d+r(r-r_0)(g(Y)-1)$.    A natural  question is \\

\centerline{\it How many points  has $Quot^E_{{\cal E}_0}$?}\vspace{3mm}

This question has a beautiful, simple  geometric interpretation: equivalently, one can
ask:\\

{\it How many holomorphic subsheaves ${\cal E}\hookrightarrow {\cal E}_0$ of
rank $r$ and maximal possible degree $d$ exist, when ${\cal E}_0$ is general?}\\

Using our computation of the gauge theoretical Gromov-Witten invariants, one gets
easily the following answer\\

\centerline{\it $Quot^E_{{\cal E}_0}$ has $r_0^{g(Y)}$ points, when $r=1$ 
and the multiplicities   are taken into account.} 
\vspace{3mm} 
Note that Lange [La] obtained earlier the inequality 
$$\#(Q)\leq 2^{g(Y)}$$
for $r_0=2$  using algebraic geometric methods, whereas Oxbury [Oxb] obtained the
equality in the smooth case  for arbitrary $r_0$.

Results in the non-abelian case  $\rk(E)>1$ were recently announced by
Lange-Newstead [LN]: \\

{\it  Suppose $g(Y)=2$, $r_0>2$ and $d$ is odd. Then $Quot^E_{{\cal E}_0}$ has
$\frac{r_0^3}{48}(r_0^2+2)$ points, if multiplicities are taken into account. }\\

Note finally that, in a recent preprint, Holla [Ho]  announced the computation of all
twisted equivariant gauge theoretical Gromov-Witten invariants of this SFAPS, which
implies  a formula for the length of any zero-dimensional Quot space $Quot^E_{{\cal
E}_0}$ of vanishing expected dimension.\\ \\
3. We mention finally another interesting result [OT2] which states that  
\\  

{\it  The the full
Seiberg-Witten invariants [OT1] of   ruled surfaces can be naturally identified with
certain gauge theoretical  Gromov-Witten invariants  of the abelian SFPAS
$\sigma_{\mu_t}=(\Hom(\C,\C^{r_0}),\alpha_{\rm can}, S^1,\mu_t)$, so they
can be deduced from Theorem 4.5.}\\

\paragraph{2. The SFPAS which yields the flag manifolds}\hfill{\break}

We come back to the SFPAS
$$(\bigoplus_{i=1}^m \Hom(V_i,V_{i+1}),  \alpha_{\rm can}, \prod_{i=1}^m
\U(V_i),\mu_t)
$$
which yields the flag manifolds with their natural symmetry.

The gauge theoretical problem associated to this SFPAS  is the following:\\

Fix $(m+1)$ Hermitian bundles $E_1,\dots, E_m, E=E_{m+1}$ over a Riemannian
surface $(Y,g)$, and fix a Hermitian connection $A_{m+1}\in{\cal A}(E)$.

Our moduli space is the space of gauge equivalence classes of systems
$$((A_1,\dots,A_m),(\varphi_1,\dots,\varphi_m))\ ,
$$
where $A_i\in{\cal A}(E_i)$, $\varphi_i\in A^0\Hom(E_i,E_{i+1})$, which solve
the equations
$$\left\{\begin{array}{cccc}
\bar\partial_{A_i,A_{i+1}}(\varphi_i)&=&0\\
i \Lambda
F_{A_i}+\frac{1}{2}(\varphi_{i-1}\circ\varphi_{i-1}^*
-\varphi_i^*\circ\varphi_i)&=&-t_i\id_{E_i},&1\leq i\leq m\ . 
\end{array}\right.
$$
Here we put of course $\varphi_{0}:=0$.  Two systems
$$((A_1,\dots,A_m),(\varphi_1,\dots,\varphi_m))\ ,\
((A'_1,\dots,A'_m),(\varphi'_1,\dots,\varphi'_m))
$$
are  equivalent if they are conguent modulo the gauge group
$${\cal G}=\prod_{i=1}^m\Aut(E_i)\ .
$$

The complex geometric classification problem which correspondes to this  gauge
theoretical problem is the following:  \\

Fix a holomorphic bundle ${\cal E}={\cal E}_{m+1}$ on a complex curve $Y$. Classify
the systems
$$(({\cal E}_1,\dots,{\cal E}_m),(\varphi_1,\dots,\varphi_m))\ ,
$$
where ${\cal E}_i$ is a    holomorphic bundle of $  {\cal
C}^\infty$-type $E_i$, and $\varphi_i\in H^0({\cal E}_i^\vee\otimes{\cal
E}_{i+1})$.  Two such systems
$$(({\cal E}_1,\dots,{\cal E}_m),(\varphi_1,\dots,\varphi_m))\ , (({\cal
E}_1',\dots,{\cal E}_m'),(\varphi_1',\dots,\varphi_m'))
$$
are  equivalent if there exist holomorphic isomorphisms $f_i:{\cal
E}_i\ra{\cal E}_i'$ such that 
$$f_{i+1}\circ \varphi_i=\varphi_i'\circ f_i\  {\rm for}\ i\in\{1,\dots,m-1\}\ ,\ 
\varphi_m=\varphi_m'\circ f_m\ .$$

The $(\mu_t,g)$-stability condition which corresponds  to this complex geometric
classification problem follows easily from the universal Kobayashi-Hitchin correspondence for
K\"ahlerian SFPAS's [LT]. This stability condition and   applications of the Kobayashi-Hitchin
correspondence in this interesting case, will be addressed  in a future article.

\paragraph{3. The SFPAS which yields the toric varieties}\hfill{\break}

Let  $Y$ be a closed connected oriented   real surface. The data
of a
$[S^1]^r$-bundle $\hat P$ on $Y$ is equivalent to the data of a system 
$L=(L_j)_j$ of $r$ Hermitian line bundles  on $Y$.  

Fix an  integer  matrix $v\in M_{m,r}(\Z)$ of rank $m$    with  the properties {\bf
P}$_1$, {\bf P}$_2$ of chapter 1, and
let $w:[S^1]^r\ra [S^1]^m$ be the associated epimorphism. Put
$$L_i^0:=\otimes_{j=1}^r [L_j^{\otimes v^i_j}]\ .
$$

The gauge theoretical problem associated with the SFAPS
$$\sigma_{\mu_\tg}=(\C^r,\alpha_{\rm can},K_w,\mu_\tg)\ ,\ \tg\in\
\coker [(v\otimes\id_\R)^*]$$ 
becomes:
\\

Choose a Riemannian metric $g$ on $Y$, a 
system of Hermitian connections $A^0=(A_i^0)_{1\leq i\leq m}$  on
$(L_i^0)_{1\leq i\leq m}$,  and put $\pg=(J_{\rm can},g,A^0)$. The moduli
space associated with  the SFAPS
$(\C^r,\alpha_{\rm can},K_w,\mu_\tg)$ and these parameters is the space    
 ${\cal M}^L_{\pg}(\sigma_{\mu_\tg})$ of equivalence classes of systems
$(A_j,\varphi_j)_j$ where 
\begin{itemize}
\item    $A_j$  is a Hermitian connection  on
$L_j$  for every $j\in\{1,\dots,r\}$, and
$$\otimes_{j=1}^r [A_j^{\otimes v^i_j}]=A_i^0\ ,\ \forall  \
i\in\{1,\dots,  m\}\ ; $$
\item  $\varphi_j$ is an $A_j$-holomorphic section    in
$L_j$ for every
$j\in\{1,\dots, r\}$\ ,  
\end{itemize}
such that $(A_j,\varphi_j)_{1\leq j\leq r}$ solves the vortex-type
equation
$$p_v\left[(i\Lambda_gF_{A_j}-2\pi\deg(L_j)+\frac{1}{2}|\varphi_j|^2)_j\right]=
\tg\ . 
$$

Two such systems are considered equivalent if they are in the same orbit with
respect to  the natural action of the gauge group  
$${\cal G}={\cal
C}^{\infty}(Y,K_w)=\{(f_1,\dots f_r)\in{\cal C}^\infty(Y,S^1)^r\ |\
\prod\limits_{j=1}^r  f_j^{v^i_j}=1\ \forall i\in\{1,\dots,m\}\} \ .
$$

The complex geometric analoga of these  notions  are the  following

\begin{dt} Let ${\cal L}^0_i$ be a holomorphic structure on $L^0_i$ for every
$i\in\{1,\dots,m\}$, and put ${\cal L}^0:=({\cal L}^0_i)_{1\leq i\leq m}$. A
holomorphic system of type
$(L,v,{\cal L}^0)$  is  a system 
$$(({\cal L}_j)_{1\leq j\leq r},(\varepsilon_i)_{1\leq i\leq m},(\varphi_j)_{1\leq
j\leq r})\ ,$$
where 
\begin{itemize}
\item ${\cal L}_j$ is a holomorphic line bundle   on $Y$,
\
\item
$\varepsilon_i: \bigotimes\limits_{j=1}^r [{\cal L}_j^{\otimes v^i_j}]\ra
{\cal L}_i^0$ is a holomorphic isomorphism for every $i\in\{1,\dots,m\}$   
 and  there exist differentiable isomorphisms  ${  L}_j\textmap{g_j}{\cal
L}_j$  with
$$\varepsilon_i\circ \bigotimes\limits_{j=1}^r [g_j^{\otimes v^i_j}]=\id_{L_i^0}
$$
\item   $\varphi_j\in H^0({\cal
L}_j)$.
\end{itemize}

An isomorphism between two such systems 
$$(({\cal L}_j)_{1\leq j\leq r},(\varepsilon_i)_{1\leq i\leq
m},(\varphi_j)_{1\leq j\leq r})\ ,\ (({\cal L}'_j)_{1\leq j\leq
r},(\varepsilon_i')_{1\leq i\leq m},(\varphi_j')_{1\leq j\leq r})$$
is  a system of holomorphic isomorphisms
$(u_j)_j$, $u_j:{\cal L}_j\ra{\cal L}'_j$ such that $\varphi'_j=u_j(\varphi_j)$ and
$\varepsilon'_i\circ [\otimes_j  u_j ^{\otimes v^i_j}]=\varepsilon'$.
\end{dt}

In order to be able to introduce the corresponding stability condition for our problem, we
need some notations: Let
$\Sigma$ be a complete simplicial fan with $\Sigma(1)\subset\{\R_{\geq
0}v^1,\dots,\R_{\geq 0}v^r\}$. For every system $T=(T_j)_{1\leq j\leq r}$   of
$r$ complex vector spaces, we put
$$U(\Sigma,T):=\{\tau\in\oplus_{j=1}^r T_j|\ 
\exists\sigma\in\Sigma \
\ {\rm such\ that}\ \tau_j\ne0\
\forall   j\in\{1,\dots,r\}\ {\rm for\ which} $$
$$ \R_{\geq 0}v^j \ {\rm
is\ not\ a\ face\ of}\ \sigma  \}\ .
$$
\begin{dt}  Let $\Sigma$ be a complete
simplicial fan in $\R^m$ such that  
$$\Sigma(1)\subset \{\R_{\geq
0}v^1,\dots,\R_{\geq 0}v^r\}\ .$$ 
A system  $(({\cal L}_j)_{1\leq j\leq r},(\varepsilon_i)_{1\leq i\leq
m},(\varphi_j)_{1\leq j\leq r})$ of type $(L,v,{\cal L}^0)$ is
$\Sigma$-stable   if one of the following equivalent conditions is
satisfied:
\begin{enumerate}
\item There exists a non-empty Zariski open  set
$Y_0\subset Y$ such that for every $y\in Y_0$ one has 
$(\varphi_1(y),\dots,\varphi_r(y))\in U(\Sigma,{\cal L}_y)$. 
\item $\varphi\in U(\Sigma,H^0({\cal L}))$, where $H^0({\cal L}):=(H^0({\cal
L}_j))_{1\leq j\leq r}$.
\end{enumerate}
\end{dt}

We denote by ${\cal M}^{\Sigma{\raisebox{0.0ex}{-}}\rm
st}_{{\cal L}_0}(L,v)$ the moduli space of $\Sigma$-stable systems of type
$(L,v,{\cal L}^0)$.
With these notations, we can now state   the Kobayashi-Hitchin correspondence for our 
problem:
\begin{thry} \ [OT3]  Let $A^0=(A_i^0)_i$,
$A^0_i\in {\cal A}(L_i^0)$ be a fixed system of  Hermitian connections and let
${\cal L}_i^0$ be the corresponding holomorphic structures. Let
$\Sigma$ be a complete simplicial fan with
$\Sigma(1)\subset\{\R_{\geq 0}v^1,\dots,\R_{\geq 0}v^r\}$, and  $\tg\in 
K_0(\Sigma)$. Then   there is
a natural isomorphism of real analytic orbifolds
$${\cal M}_{\pg}^L(\sigma_{\mu_\tg})\simeq {\cal
M}^{\Sigma{\raisebox{0.0ex}{-}}\rm st}_{{\cal L}^0}(L,v)\ .
$$
\end{thry}

Using this theorem  one  obtains the following important results [OT3]:
\begin{thry} Suppose that $-i\tg$ is a regular value of the standard moment map
$\mu_0$.  Then
\begin{enumerate}
\item (complex geometric interpretation) The moduli space ${\cal
M}_{\pg}^L(\sigma_{\mu_\tg})$ is a toric fibration over an abelian  variety $P$ of
dimension $g(Y)(r-m)$.
\item (embedding theorem) The moduli space ${\cal
M}_{\pg}^L(\sigma_{\mu_\tg})$ can be identified with  the vanishing locus of a  
section $\sigma$ in a split holomorphic   bundle  
${\cal E}$ over the total space of a   locally trivial holomorphic toric fibre
bundle
$T$ over $P$.
\end{enumerate}
\end{thry}

In the special case $Y=\P^1$, ${\cal L}_i^0={\cal O}_{\P^1}$, the complex
geometric interpretation was previously obtained and used in [W2], [MP]. 

Using the "embedding theorem", one can endow the moduli space ${\cal
M}_{\pg}^L(\sigma_{\mu_\tg})$ with a natural
{\it algebraic geometric} virtual fundamental class, namely the localized Euler class
$[Z(\sigma)]$ of the bundle ${\cal E}$ [F]. Roughly speaking, this class -- which is an element
in the Chow group $A_{\dim(T)-\rk({\cal E})}(Z(\sigma))$ -- is obtained by "intersecting" the
image of
$\sigma$ with the zero section of
${\cal E}$ in the (smooth!) total space of this vector bundle.

\begin{thry} (comparison theorem)  The image of the algebraic geometric virtual
fundamental class $[Z(\sigma)]$ of
$Z(\sigma)$ under the cycle map 
$$cl:A_{*}(Z(\sigma))\map H_{*}(Z(\sigma))$$  
agrees with the  gauge theoretical virtual fundamental class $[{\cal
M}_{\pg}^L(\sigma_{\mu_\tg})]^{\rm vir}$ via the embedding theorem above.
\end{thry}

This (rather difficult) result is the first explicit verification of the conjecture stated in
section 3.

The computation of the non-twisted gauge theoretical Gromov-Witten invariants in
the case $g(Y)=0$ was  obtained in [CS].  Results in the "higher genus case" can
be found in [Ha], where the comparison Theorem 4.10 plays an important role.

\paragraph{4. The   SFPAS of Str{\o}mme
triples}\hfill{\break}

Now we consider again the SFPAS
$$
\sigma_{\mu_{s,t}}=(\Hom(U,V)^{\oplus 2}\oplus\Hom(W,V),   \alpha_{\rm
can},\U(U)\times\U(V),\mu_{s,t})\ ,
$$
where
$$\mu_{s,t}(k,l,m)=\frac{i}{2}(k^*\circ k+l^*\circ l ,-k\circ k^*-       
l\circ l^*-m\circ m^*)+i(-t\id_{U}, s\id_{V})\ .
$$

The gauge theoretical problem associated to this SFPAS is the following:\\

Let $E$, $F$, $H_0$ be fixed differentiable Hermitian  bundles  on $Y$, $A_0$
  a fixed Hermitian connection on $H_0$, and $g$ a Riemannian metric on $Y$.
We are interested in the moduli space ${\cal M}_{A_0,g}(E,F,H_0)$ of
equivalence classes of systems
$(A,B,(k,l,m))$, where 
\begin{itemize}
\item $A\in{\cal A}(E)$, $B\in{\cal A}(F)$, 
\item $k$, $l\in A^0\Hom(E,F)$, $m\in A^0\Hom(H,F)$,
\end{itemize}
such that
\begin{itemize}
\item $k$, $l$ are $\bar\partial_{A,B}$-holomorphic, $m$ is
$\bar\partial_{A_0,B}$-holomorphic,
\item the equations
$$\left\{
\begin{array}{ccc}
i\Lambda F_A-\frac{1}{2}(k^*\circ k+l^*\circ l)&=&-t\id_E\\
i\Lambda F_B+\frac{1}{2}(k\circ k^*+l\circ l^*+m\circ m^*)&=&s\id_E
\end{array}
\right.
$$
\end{itemize}
are satisfied.

Two such systems are considered equivalent if they are congruent modulo the action of the
gauge group
$${\cal G}=\Aut(E)\times \Aut(F)\ .
$$

The corresponding complex geometric classification problem   is the following:\\

Let ${\cal H}_0$ be the holomorphic structure on $H_0$ defined by
$A_0$. Classify  systems $({\cal E},{\cal F},(k,l,m))$ where 
\begin{itemize}
\item ${\cal E}$, ${\cal F}$ are
holomorphic bundles which are differentiably isomorphic with $E$ and $F$
respectively,
\item 
$(k,l)\in H^0({\cal E}^\vee\otimes{\cal F})^{\oplus 2}$,   $m\in H^0({\cal
H}_0^\vee\otimes{\cal F})$.
\end{itemize}

Such a system will be called  holomorphic systems of type $(E,F,{\cal H}_0)$. Two
holomorphic systems 
$({\cal E},{\cal F},(k,l,m))$,
$({\cal E}',{\cal F}',(k,l,m))$ of type $(E,F,{\cal H}_0)$ are   equivalent  if
there exist holomorphic isomorphisms
$f:{\cal E}\ra{\cal E}'$, $g:{\cal F}\ra{\cal F}'$ such that $k'\circ f=g\circ k$,
$l'\circ f=g\circ l$,
$g\circ m=m'$.

\begin{dt} A   holomorphic system   $({\cal E},{\cal
F},(k,l,m))$ of type $(E,F,{\cal H}_0)$ is
$(s,t)$-stable if and only if for every  pair of   subsheaves ${\cal E}_1\subset  {\cal
E}$, ${\cal F}_1\subset {\cal F}$ with torsion-free  quotients ${\cal
E}^2:={\cal E}/{\cal E}_1$, ${\cal F}^2 :={\cal F}/{\cal
F}_1$ one has
\begin{enumerate}
\item When $({\cal E}_1,{\cal F}_1)\ne(0,0)$ and
$$k({\cal E}_1)+l({\cal
E}_1)\subset {\cal F}_1\ ,
$$
the following inequality holds
$$\frac{2\pi}{(n-1)!}\left[-\deg({\cal E}_1)-\deg({\cal F}_1) \right]+
 \Vol_g(X)[s\rk({\cal F}_1)-t\rk({\cal E}_1)]> 0\ .$$
\item When $({\cal E}^2,{\cal F}^2)\ne (0,0)$ and 
$$k({\cal E}_1)+l({\cal E}_1)\subset {\cal F}_1\ ,\ \im(m)\subset {\cal F}_1\ ,
$$
the following inequality holds
$$\frac{2\pi}{(n-1)!}\left[\deg({\cal E}^2)+\deg({\cal F}^2)
\right]+\Vol_g(X)[t\rk({\cal E}^2)-s\rk({\cal F}^2)]> 0\ .$$
\end{enumerate}
\end{dt}
There is a moduli space ${\cal M}^{(s,t)-{\rm st}}_{{\cal H}_0}(E,F)$ classifying
stable  systems of a fixed type $(E,F,{\cal H}_0)$. In this case the Kobayashi-Hitchin 
correspondence for our moduli problem has the following form:
\begin{thry} \ [LOT]  There is a
natural isomorphism of real analytic spaces
$${\cal M}_{A_0,g}(E,F,H_0)^*\simeq {\cal M}^{(s,t)-{\rm st}}_{{\cal H}_0}(E,F)\ .
$$
\end{thry}

\paragraph{5.  SFPAS's associated with quiver classification problems}\hfill{\break}

 The moduli problem defined by  the standard twisted quiver factorization problem
(see chapter 1, paragraph 5) and the associated  Kobayashi-Hitchin correspondence was
studied  in [AlPr1], [AlPr2] for 
arbitrary  K\"ahler manifolds.  We recall that this case
corresponds to the SFPAS of the form
$$(\bigoplus_{a\in A}\Hom(W_{s(a)},W_{t(a)}\otimes W_a^0), \alpha_{\rho,r},\prod_{v\in
V} \U(W_v), \mu_\tg)\ ,
$$
where $Q=(V,A,s,t)$ is a quiver, $W_v$ and $W_a^0$ are Hermitian vector
spaces indexed by $v\in V$ and $a\in A$,   $\alpha_{\rho,r}$ is the canonical
representation of 
$$\hat K:=\prod_{v\in V} \U(W_v)\times \prod_{a\in A}\U(W_a^0)
$$
on $F=\bigoplus_{a\in A}\Hom(W_{s(a)},W_{t(a)}\otimes W_a^0)$, and $\mu_\tg$ is the
moment map for the $K$-action on $F$.
  
The  moduli problem corresponding to this SFPAS is the following:

Fix Hermitian bundles $(E_v)_{v\in V}$, $(E_a^0)_{a\in A}$ on a Riemann surface
$(Y,g)$, and fix  connections $A^0_a$ on $E_a^0$.  Classify all systems
$((A_v)_{v\in V},(\varphi_a)_{a\in A}))$,
where 
\begin{itemize}
\item $A_v\in{\cal A}(E_v)$ is a Hermitian connection,
\item $\varphi_a\in
A^0(\Hom(E_{s(a)},E_{t(a)}\otimes E_a^0))$  
\end{itemize}
such that
$$\left\{\begin{array}{ccc}
\bar\partial_{A_{s(a)}, A_{t(a)}\otimes A_a^0}\ (\varphi_a)&=&0\ ,\ \forall a\in A,\\
i\Lambda F_{A_v}+\frac{1}{2}(\sum\limits_{t(a)=v}\tr_{E_a^0}   
 (\varphi_a\circ\varphi_a^*)-\sum\limits_{s(a)=v} 
 (\varphi_a^*\circ\varphi_a))&=&-t_v\id_{E_v}\ .
\end{array}
\right.
$$

For the corresponding complex geometric classification problem and the Kobayashi-Hitchin
correspondence in this case we refer the reader to [AlPr2].  

The  SFAPS associated with    {\it general}
(twisted) quiver factorization problems is  more difficult and leads to  interesting    moduli
problems.  The corresponding complex geometric classification problem and the stability
condition follows from the universal Kobayashi-Hitchin correspondence for K\"ahlerian   
SFPAS's [LT].

 {\small
Authors addresses: \vspace{2mm}\\
Ch. Okonek,   Institut f\"ur Mathematik, Universit\"at Z\"urich,  
Winterthurerstrasse 190, CH-8057 Z\"urich, Switzerland,  e-mail:
okonek@math.unizh.ch\vspace{1mm}\\
 A. Teleman, LATP, CMI,   Universit\'e de Provence,  39  Rue F. J.
Curie, 13453 Marseille Cedex 13, France,  e-mail: teleman@cmi.univ-mrs.fr 
    ,  and\\ Faculty of Mathematics, University of Bucharest, Bucharest,
Romania }

\end{document}